\documentclass[12pt]{extarticle} 
\usepackage[top=1in, bottom=1in,left=1in,right=1in]{geometry}

\usepackage{makeidx, amsfonts, amsthm, amssymb, amsmath, graphicx, subfig, bbm, color,epstopdf, wrapfig, caption, listings, float, mathtools, hyperref, url, dsfont}

\usepackage[capitalize,nameinlink]{cleveref}

\allowdisplaybreaks

\usepackage[shortlabels]{enumitem}
\setlist{leftmargin=*,itemsep=0.25\itemsep,parsep=0.35\parsep,topsep=0.25\topsep,partopsep=0.29\partopsep}

\newtheorem{theorem}{Theorem}
\newtheorem{other}{Theorem}
\newtheorem{definition}[other]{Definition}
\newtheorem{lemma}[other]{Lemma}

\newtheorem*{remark*}{Remarks}
\newtheorem{conjecture}[theorem]{Conjecture}

\renewenvironment{proof}[1][\proofname]{ { \it\bfseries #1: }}{\qed}

\newcommand{\cN}{{\mathcal{N}}}

\newcommand{\field}[1]{\mathbb{#1}}
\newcommand{\R}{\field{R}}

\newcommand{\p}{\field{P}}
\DeclareMathOperator{\tr}{tr}

\DeclareMathOperator{\E}{\mathbb{E}}

\newcommand{\mat}[2][rrrrrrrrrrrrrrrrrrrrrrrrrrrrrrrr]{\left[ \begin{array}{#1} #2 \\ \end{array}\right]}

\numberwithin{equation}{section}

\providecommand{\keywords}[1]
{
  \small	
  \textbf{\textit{Keywords: }} #1
}

\title{Invariance principle of random matrix for the norm}
\author{Juntao Duan\footnote{corresponding author, School of Mathematics, Georgia Institute of Technology, 
    \href{mailto:juntaoduan@gmail.com}{juntaoduan@gmail.com}}
    ,\quad  
    Ionel Popescu\footnote{ University of Bucharest, Faculty of Mathematics and Computer Science, Institute of Mathematics of the Romanian Academy, \href{mailto:ioionel@gmail.com }{ioionel@gmail.com}}
    ,\quad  
    Heinrich Matzinger\footnote{School of Mathematics, Georgia
 Institute of Technology, \href{mailto:matzi@math.gatech.edu}{matzi@math.gatech.edu}}
    }

\begin{document}

\large
\maketitle

\begin{abstract}
	 Johnson-Lindenstrauss guarantees certain topological structure is preserved under random projections when project high dimensional deterministic vectors to low dimensional vectors. In this work, we try to understand how random matrix affect norms of random vectors. In particular we prove the distribution of the norm of random vector $X \in \mathbb{R}^n$, whose entries are i.i.d. random variables, is preserved by random projection $S:\mathbb{R}^n \to \mathbb{R}^m$. More precisely, 
\[
\frac{X^TS^TSX - mn}{\sqrt{\sigma^2 m^2n+2mn^2}} \xrightarrow[\quad m/n\to 0 \quad ]{ m,n\to \infty } \mathcal{N}(0,1)
\]
We also prove a concentration of the random norm transformed by either random projection or random embedding. Overall, our results showed random matrix has low distortion for the norm of random vectors with i.i.d. entries.
\end{abstract}

\section{Introduction}

     Due to the internet boom and computer technology advancement in the last few decades, data collection and storage have been growing exponentially. With 'gold' mining demand on the enormous amount of data reaches to a new level, we are facing many technical challenges in understanding the information we have collected. In many different cases, including text and images, data can be represented as points or vectors in high dimensional space. On one hand, it is very easy to collect more and more information about the object so that the dimensionality grows quickly. On the other hand it is very difficult to analyze and create useful models for high dimensional data due to several reasons including computational difficulty as a result of curse of dimensionality and high noise to signal ratio. It is therefore necessary to reduce the dimensionality of the data while preserving the relevant structures. 
     
     The celebrated Johnson-Lindenstrauss lemma \cite{johnson1984extensions} states that random projections can be used as a general dimension reduction technique to embed topological structures in high dimensional Euclidean space into a low dimensional space without distorting its topology. 
     Since then random projections has been found very useful in many applications such as signal processing and machine learning.  For example fast Johnson-Lindenstrauss random projections is used to approximate K-nearest neighbors to speed up computation \cite{indyk1998approximate, ailon2009fast}.  Random sketching uses random projection to reduce sample sizes in regression model and low rank matrix approximation \cite{woodruff2014sketching}. Random projected features can be used to create low dimensional base classifiers which are combined as robust ensemble model \cite{cannings2017random}. Practitioners found applications of random projection in privacy and security \cite{liu2005random}.
Let us first recall the Johnson-Lindenstrauss lemma \cite{burr2018optimal}.
  \begin{lemma}[Johnson and Lindenstrauss]\label{lemma:Johnson-Lindenstrauss lemma}
       Given a set of vectors $\{u_1, \cdots, u_k\}$ in $\R^n$, for any $m \ge 8 \varepsilon^{-2} \log k$, there exists a linear map $A:\R^n \to \R^m$ such that 
       \[
       (1-\varepsilon)\|u_i -u_j\| \le 
       \|Au_i -Au_j\| \le (1+\varepsilon)\|u_i -u_j\|
       \]
  \end{lemma}
  
  Given two fixed vectors $X_1, X_2 \in \R^n$,  by Johnson-Lindenstrauss lemma,  we can find a random projections $A: \R^n \to \R^m$ such that the projected distance $\|AX_1-AX_2 \|$ has only a small distortion of the original distance $\|X_1-X_2 \|$. More precisely, 
\begin{align}\label{eqn:J-L lemma norm error order}
    \left[1- O(\frac{1}{\sqrt{m}}) \right] \|X_1-X_2\|^2 \le \| A(X_1-X_2)\|^2 \le \left[1+ O(\frac{1}{\sqrt{m}}) \right] \|X_1-X_2\|^2
\end{align}
Equivalently, this property can be reformulated as random projections preserves the inner product of two vectors (equivalence can be obtained by elementary computation and polarization identity). 
Namely given $X_1, X_2$ two vectors in the unit ball of $\R^n$ ($\|X_1\|\le 1, \|X_2\|\le 1$) , then there is a random projection $A:\R^n \to \R^m$ such that 
\begin{align}\label{eqn:J-L lemma inner product}
    | \langle AX_1, AX_2 \rangle - \langle X_1, X_2 \rangle | \le O(\frac{1}{\sqrt{m}})
\end{align}
For general vectors not in the unit ball, the bound on the right hand side has the norms as a factor \[
| \langle AX_1, AX_2 \rangle - \langle X_1, X_2 \rangle | \le O(\frac{1}{\sqrt{m}})\|X_1\| \|X_2\|
\]

The natural extension is to consider random vectors $X_1, X_2$. The question becomes how random projections affect random vectors. Suppose $X_1$ and $X_2$ are independent. After applying a random projection or embedding, independent random vectors become strongly dependent. Does this mean random projection is inferior to be used for dimension reduction of random vectors?  \cite{duan2021invariance} showed there is an invariance phenomenon, namely the distribution of inner product of two independent with i.i.d. entries is preserved by random matrix. 


  In this work we try to focus on inner product of two dependent random vectors. In particular we only address the case that the two vectors are the same. Therefore we will obtain an invariance of distribution of the randomly projected norm. To put this in a high level perspective, we shall introduce the full inner product structure.


\subsection{The full inner product structure}\label{sec:RP full structure}

Suppose there are $p$ independent random vectors $X_1, \cdots X_p$ with i.i.d. entries (say $x$), and a random projection matrix $S$ of dimension $m\times n$ with i.i.d. entries (say $s$). The full inner product structure of those vectors and $S$-projected vectors concerns two random matrices each collecting the following inner products:
 \[
 \langle X_i, X_j \rangle, \;  \langle SX_i, SX_j \rangle \quad \forall 1\le i, j \le p
 \]
 
 In the previous work (\cite{duan2021invariance}), the cases that $i\neq j$ are addressed.
 We will focus on the remaining $p$ diagonal terms in this work.  Namely, we will try to understand how random projection affects the distribution of the norm of a random vector $X$,
 \[
 \langle X, X \rangle, \;  \langle SX, SX \rangle 
 \]


 Before getting into technicality, we will first use techniques similar to the proof of Johnson-Lindenstrauss lemma to obtain a Bernstein type concentration result (\cref{thm:RP norm concentration}) under sub-Gaussian assumptions in \cref{sec:concen norm by JL}.
Such concentration properties allow one to analyze the problem from error control perspective. In particular, this shows the randomly projected norm admits a sub-Gaussian behavior.

In \cref{sec:RP  distribution of norm}, we would go one step further to deal with the distribution directly and we show the distribution of norm is invariant under random projections (\cref{thm:rand_proj norm CLT}), namely the random projected norm converges to normal distribution after properly centered and scaled. 
\[
\frac{\langle SX, SX \rangle  - mn}{\sqrt{\sigma^2 m^2n+2mn^2}} \xrightarrow[ m/n\to 0 ]{ m,n\to \infty } \cN(0,1), \quad   \text{ where } \sigma^2 = \E x^4 -1 
\]
Lastly, we outline open questions and possible future work in \cref{sec:Open questions}.
 

\section{Concentration of projected or embedded norm for sub-Gaussian variables}\label{sec:concen norm by JL}
 
The purpose of this section is to show the randomly projected norm is concentrated around both the original random norm and expectation of the norm. Concentration inequalities concern the tails of a random quantity deviates from its mean, which are very powerful tools in many applications (\cite{ledoux2001concentration, boucheron2013concentration, vershynin2018high}). Johnson-Lindenstrauss  \cref{lemma:Johnson-Lindenstrauss lemma} itself is a result of concentration inequality for sub-Gaussian random matrix over fixed vectors (see \cite{boucheron2013concentration, vershynin2018high}). Concentration properties of random quadratic forms involving either deterministic vectors or deterministic matrix have been well-studied in the literature (see \cite{tao2011topics, rudelson2013hanson, tropp2015introduction, vershynin2018high}). Most of the existing results control the tail probability of the distortion by the matrix or expected distortion quantitatively. We shall use similar techniques to prove the concentration of randomly projected norm of sub-Gaussian random vectors.  First let us recall some properties of sub-Gaussian  random variables.

\begin{definition}
We say $X$  is a sub-Gaussian random variable if there is $v>0$ such that  
\[
\E e^{\lambda X} \le e^{\frac{\lambda^2v}{2}}
\]
\end{definition}
Using Markov inequality, we can easily see sub-Gaussian random variable admits the tail probability
$$ \p(|X|>t)\le 2 e^{-\frac{t^2}{2v}}, \quad \forall t \ge 0 $$

For now let us assume sub-Gaussian random variable $X$ is centered and standardized, namely $\E X=0, v=1$. It is not hard to verify sub-Gaussian tail property implies moments bounds (see section 2.3 of \cite{ boucheron2013concentration}). 
\begin{align*}
    \E X^{2q} & =\int_0^{\infty} \p(X^{2q}>s) ds \\
    & = \int_0^{\infty} q t^{q-1}\p(X^{2}>t) dt \\
    &\le \int_0^{\infty} q t^{q-1} 2 e^{-t/2} dt \\
    &\le 2^{q+1} q!
\end{align*}
This will allow us to compute moment generating function of $X^2$, and some useful bounds to be used later. Firstly,
\begin{align}\label{eqn:sub-gaussian square exp-moment bound noncenter}
    \E e^{\lambda X^2} 
   & = 1+ \sum_{q=1}^{\infty} \frac{\lambda^{q} \E X^{2q} }{q!} \nonumber \\
   & \le 1+ \sum_{q=1}^{\infty} \lambda^{q} 2^{q+1} \nonumber \\
   & = \frac{1+2\lambda}{1-2\lambda}, \qquad  \forall \; \lambda< \frac{1}{2} \nonumber \\
   &\le e^{5\lambda}, \qquad  \forall \; \lambda< \frac{1}{5}
\end{align}

Secondly, let $X'$ be an independent copy of $X$.
\begin{align*}
    \E e^{\lambda (X^2 -1)} \E e^{-\lambda (X'^2 -1)} 
   & =  \E e^{\lambda [(X^2 -1)-(X'^2-1)} 
   = 1+\sum_{q=1}^{\infty} \frac{\lambda^{2q} \E (X^2-X'^2)^{2q}}{(2q)!}
\end{align*}
Notice  by Minkowski's inequality we have $\E (X^2-X'^2)^{2q}\le 2^{2q}\E X^{4q}$, and  by Jensen's inequality $ \E e^{-\lambda (X'^2 -1)} \ge e^{-(\E X'^2-1)}=1$. Therefore
\begin{align}\label{eqn:sub-gaussian square exp-moment bound}
    \E e^{\lambda (X^2 -1)} 
   & \le 1+ \sum_{q=1}^{\infty} \frac{\lambda^{2q} 2^{2q} 2^{2q+1} (2q)!}{(2q)!} \nonumber \\
   & \le \frac{1+16\lambda^2}{1-16\lambda^2} , \qquad \forall \; \lambda< \frac{1}{4}\nonumber \\
   & \le e^{40\lambda^2}, \qquad \forall \; \lambda< \frac{1}{5}
\end{align}

From a high level, these properties (\cref{eqn:sub-gaussian square exp-moment bound noncenter}, \cref{eqn:sub-gaussian square exp-moment bound}) of $X^2$  are expected since it is a sub-exponential random variable. For a centered sub-Gaussian $X$ with variance $1$, it is easy to see $X^2$ has sub-exponential tail decay. Because $\p(X^2> t)= \p( |X|>\sqrt{t}) \le 2 e^{-\frac{t}{2}}$. The centered version $X^2-1$ has tails shifted by a constant $1$ thus again admits sub-exponential decay.  For extensive detailed discussions and proofs of the properties, we refer to  \cite{boucheron2013concentration, vershynin2018high}.



\bigskip
\begin{theorem}\label{thm:RP norm concentration}
Given $X= [x_1,\cdots,x_n]^T$. Let $A=\frac{1}{\sqrt{m}} S$ be a $m\times n$ random matrix. Let all random variables $x_i, S_{i,j}$ are independent identically distributed with mean zero and variance one. Suppose all random variables are sub-Gaussian, then we have 
\begin{enumerate}
\item 
    \begin{equation}\label{eqn:two side tail sum squares of subgaussian}
    \p \left( \frac{1}{\sqrt{n}} \left|  \|X\|^2 -n \right|> t \right) < 2\exp\left\{ {- \min(\frac{t^2}{160}, \frac{t\sqrt{n}}{10}) }  \right\}
    \end{equation}
\item 
    \begin{equation} \label{eqn:tail of random norm to random projected norm}
     \p \left(  \frac{1}{\sqrt{n}} \left| \|AX\|^2 - \|X\|^2 \right|> \frac{\|X\|^2}{n} t \right) 
     < 2 \exp \left\{ -\min \left(   \frac{t^2 \; C}{160} , \; { \frac{t \; C \sqrt{n}}{10} }  \right) \right\}
    \end{equation}

\end{enumerate}
    
where $C= \frac{m}{n}$.

\end{theorem}

\begin{proof}

\begin{enumerate}

\item 
Notice $x_i$ are i.i.d. sub-Gaussian with variance $v=1$, then $ (x_i^2 -1)$ are  i.i.d. sub-exponential random variables. We can use Chernoff type argument (or apply Bernstein's  concentration inequality directly)  to calculate
\begin{align*}
    \p \left( \frac{1}{\sqrt{n}} \left(  \|X\|^2-n \right)> t \right) 
    & = \p \left( e^{\frac{\lambda}{\sqrt{n}} \left(  \|X\|^2 -n \right)} > e^{\lambda t} \right)  \\
    & \le e^{-\lambda t} \E \left[e^{\frac{\lambda}{\sqrt{n}} \left( \|X\|^2 -n \right)} \right] \\
    & =  e^{-\lambda t} \prod_{i=1}^n \left[ \E e^{\frac{\lambda}{\sqrt{n}} \left( x_i^2 -1 \right)} \right] 
\end{align*}
which holds for any $\lambda \ge 0$. We know by \cref{eqn:sub-gaussian square exp-moment bound}, for any $\lambda/\sqrt{n} \le \frac{1}{5}$, we have $\E e^{\frac{\lambda}{\sqrt{n}} \left( x_i^2-1 \right)} \le e^{ 40 \lambda^2/n  }$

then we obtain
\begin{align*}
    \p \left( \frac{1}{\sqrt{n}} \left( \|X\|^2-n \right)> t \right) 
    & \le  e^{-\lambda t + 40 \lambda^2} , \qquad   {\lambda} < \frac{{\sqrt{n}}}{5}
\end{align*}
Minimizing the right hand under the constraint  $\lambda \le \frac{\sqrt{n}}{5}$, we  find optimal $\lambda^*= \min(\frac{\sqrt{n}}{5}, \frac{t}{80})$.
Therefore we obtain
\begin{align}
    \p \left( \frac{1}{\sqrt{n}} \left(  \|X\|^2-n \right)> t \right) 
    & \le 
\begin{cases*}
\exp\{ {- \frac{t^2}{160} }  \},  & if $t < 16\sqrt{n}$  \\
\exp\{ { \frac{8}{5}n-\frac{t\sqrt{n}}{5} }  \} \le \exp\{ { -\frac{t\sqrt{n}}{10} } \},   & if $t\ge 16\sqrt{n}$
\end{cases*} \nonumber\\
& = \exp\left\{ {- \min(\frac{t^2}{160}, \frac{t\sqrt{n}}{10}) }  \right\}  \label{eqn:bounding one side tail of sum of squares of subgaussian}
\end{align}
Repeat the same argument for $\frac{1}{\sqrt{n}} \left( n- \|X\|^2 \right)$, we find the other half admits the same tail bound, thus we obtain \cref{eqn:tail of random norm to random projected norm} 
\begin{equation*}
    \p \left( \frac{1}{\sqrt{n}} \left| \|X\|^2-n \right|> t \right) < 2\exp\left\{ {- \min(\frac{t^2}{160}, \frac{t\sqrt{n}}{10}) }  \right\}
\end{equation*}

We see as $n\to \infty$, the tail has a Gaussian behavior namely $e^{-ct^2}$ which coincide with the CLT of $\frac{1}{\sqrt{n}} \left( n - \|X\|^2 \right) \to \cN(0,\E(x_1^2-1)^2)$. 

\item 
Now we want to quantify how much $\|AX\|^2$ deviates from $\|X\|^2$. 
\begin{align*}
    \|AX\|^2 &= \sum_{i=1}^m (A_{i,.}X)^2
    = \sum_{i=1}^m \left(\sum_{j=1}^n A_{i,j}x_j \right)^2 
    = \sum_{i=1}^m \sum_{j=1}^n \sum_{k=1}^n A_{i,j} A_{i,k}x_jx_k 
\end{align*}
We can use conditioning on $X$ and only deal with the conditional probability.
\begin{equation}\label{eqn:concen split, conditioning on X}
    \p \left( \frac{1}{\sqrt{n}} \left| \|AX\|^2 - \|X\|^2 \right|> t \right) = \E \left[ \p \left( \frac{1}{\sqrt{n}} \left| \|AX\|^2 - \|X\|^2 \right|> t \middle| X \right)  \right]
\end{equation}
In which case, we can think of $X$ being fixed when computing conditional probability. To simplify notation, define 
 \[
 y_i := \left[ \left(\sum_{j=1}^n A_{i,j}x_j \right) \middle| X \right] ,  \quad 1\le i \le m
 \]
 Therefore $\|AX\|^2 |X = \sum_{i=1}^m y_i^2$. Notice $\E y_i =0$ since $\E A_{i,j}=0$.
Moreover $A_{i,j}, A_{i,k}$ are independent if $j\ne k$, we find
\begin{align*}
    \E y_i^2 
& = \E \left[ \left(\sum_{j=1}^n A_{i,j}x_j \right)^2 \middle| X \right] \\
& = \E \left[ \sum_{j=1}^n\sum_{k=1}^n A_{i,j} A_{i,k} x_j x_k  \middle| X \right] \\
& = \sum_{j=1}^n x_j^2\E A_{i,j}^2  \\
&=  \sum_{j=1}^n x_j^2 \frac{1}{m}\E S_{i,j}^2 \\
&= \frac{1}{m} \|X\|^2
\end{align*}
And this also shows conditional expectation of projected norm is the original norm
\begin{align*}
    \E \left[\|AX\|^2  \middle| X \right]
    &= \sum_{i=1}^m \E [y_i^2|X]  
    = \|X \|^2 
\end{align*}
We may rewrite the tail of norms in terms of random variable $y_i$,
\begin{align*}
    \p \left( \frac{1}{\sqrt{n}} \left| \|AX\|^2 - \|X\|^2 \right|> t \right) = \E \left[ \p \left( \frac{1}{\sqrt{n}} \left|  \sum_{i=1}^m (y_i^2 - \E y_i^2)\right|> t  \right)  \right]
\end{align*}
 
In fact, linear combination of sub-Gaussian is still sub-Gaussian. We shall prove $y_i$ is  sub-Gaussian random variable so that we can obtain Bernstein's type inequality again by applying \cref{eqn:sub-gaussian square exp-moment bound}. 
 \begin{align*}
   \E[e^{y_i t}] 
    & =  \E e^{ t\sum_{j=1}^n A_{i,j}x_j}    = \prod_{j=1}^n \E e^{ x_j t A_{i,j}}
    = \prod_{j=1}^n \E e^{ x_j t \frac{S_{i,j}}{\sqrt{m}} } \\
    & \le \prod_{j=1}^n e^{ {x_j^2 t^2}/{2m} } \\
  & = e^{(\sum_{j=1}^n  x_j^2 t^2/2m) } \\
  &  = e^{(\|X\|^2/m)\frac{t^2}{2}}
 \end{align*}
 Then the tail probability
 \begin{align*}
     \p(y_i>t) & = \p(e^{y_i \lambda}> e^{t\lambda}) \\
     & \le e^{-t\lambda} \E e^{y_i \lambda} \\
     & \le \exp \left(-t\lambda + (\|X\|^2/m)\frac{\lambda^2}{2} \right) \quad \forall \lambda >0
 \end{align*}
 If we minimize on the right side over $\lambda>0$, we should take $\lambda =t m/\|X\|^2$. Therefore we obtain
 \[
 \p(y_i>t) \le \exp \left(- \frac{t^2}{2(\|X\|^2/m)} \right)
 \]
 Repeat the same argument for $\p(y_i<-s)$, we will obtain two-sided sub-Gaussian tail bound.
 This shows $y_i$ is sub-Gaussian with variance $v=\|X\|^2/m $. Therefore $y_i^2-\E y_i^2$ is sub-exponential.
 Then with Chernoff's method, we calculate 
 \begin{align*}
    \p \left( \frac{1}{\sqrt{n}} \left|  \sum_{i=1}^m (y_i^2 - \E y_i^2)\right|> t  \right) 
    & = \p \left( \frac{1}{\sqrt{n}} \left|  \sum_{i=1}^m (y_i^2/v - \E y_i^2/v)\right|> t/v  \right) \\
    & \le 2 e^{-\lambda t/v} \prod_{i=1}^m \E e^{\frac{\lambda}{\sqrt{n}}(\frac{y_i^2}{v} - \E \frac{y_i^2}{v})}
 \end{align*}
which holds for any $\lambda \ge 0$. Notice $y_i/v$ are centered and standardized independent sub-Gaussian random variables. By \cref{eqn:sub-gaussian square exp-moment bound}, we know for any $\lambda \le \frac{\sqrt{n}}{5}$, we have 
\[
\E e^{\frac{\lambda}{\sqrt{n}}(\frac{y_i^2}{v} - \E \frac{y_i^2}{v})} \le   e^{ 40 \lambda^2/n }
\]
then we obtain
 \begin{align*}
    \p \left( \frac{1}{\sqrt{n}} \left|  \sum_{i=1}^m (y_i^2 - \E y_i^2)\right|> t  \right) 
    & \le 2  \exp \left\{ -\lambda t\frac{m}{\|X\|^2} + 40\lambda^2 \frac{m}{n}  \right\} \quad \forall \lambda \le \frac{\sqrt{n}}{5}
 \end{align*}
Optimize the right hand  side under the constraint  $\lambda \le \frac{\sqrt{n}}{5}$, we  find optimal $\lambda^*= \min(\frac{\sqrt{n}}{5}, \frac{t \;n}{80\|X\|^2})$.
Therefore we obtain
\begin{align*}
\p \left( \frac{1}{\sqrt{n}} \left|  \sum_{i=1}^m (y_i^2 - \E y_i^2)\right|> t  \right) & \le 
    \begin{cases*}
    2\exp\{ {- \frac{t^2 \; mn}{160(\|X\|^2)^2} }  \},  & if $t < 16\|X\|^2/\sqrt{n}$  \\
    2\exp\{ { \frac{8m}{5}-\frac{t\sqrt{n} \; m}{5 \|X\|^2} }  \} ,   & if $t\ge 16\|X\|^2/\sqrt{n}$
    \end{cases*}  \\ 
    & \le 
    \begin{cases*}
    2\exp\{ {- \frac{t^2 \; Cn^2}{160(\|X\|^2)^2} }  \},  & if $t < 16\|X\|^2/\sqrt{n}$  \\
    2\exp\{ -\frac{t\sqrt{n} \; Cn}{10 \|X\|^2}   \} ,   & if $t\ge 16\|X\|^2/\sqrt{n}$
    \end{cases*}  
\end{align*}
where $C:=\frac{m}{n} \ge 0$ is a constant.
Therefore 
\begin{align*} 
\p \left( \frac{1}{\sqrt{n}} \left|  \sum_{i=1}^m (y_i^2 - \E y_i^2)\right|> t  \right) & \le 2\exp \left\{ -\min \left(   \frac{t^2 \; Cn^2}{160(\|X\|^2)^2} , \; { \frac{t\sqrt{n} \; Cn}{10 \|X\|^2} }  \right) \right\}
\end{align*}
Replacing $t$ by $\|X\|^2t/n$, we obtain 
\begin{equation*}
   \p \left(  \frac{1}{\sqrt{n}}\left|  \sum_{i=1}^m (y_i^2 - \E y_i^2)\right|> \frac{\|X\|^2 }{n} t  \right) 
   \le 2 \exp \left\{ -\min \left(   \frac{t^2 \; C}{160} , \; { \frac{t \sqrt{n}\; C}{10} }  \right) \right\}
\end{equation*}

Taking expectation with respect to $X$ we obtain \cref{eqn:tail of random norm to random projected norm}.

\end{enumerate}
\end{proof}

\begin{remark*}
The first tail bound implies $\|X\|^2$ is close to $n$. The second tail bound implies $\|AX\|^2$ is close to $\|X\|^2$ and thus also close to $n$. Later next section we will obtain a CLT type result for
\[
   \|AX\|^2  =\frac{1}{m} X^TSSX 
\]
  This concentration result actually suggests the centered and rescaled projected norm has a tail that decays at a sub-Gaussian rate when $t\le \sqrt{n}$ provided the random variables are originally sub-Gaussian. Thus a Gaussian limit (though without assuming random variables are sub-Gaussian) which we will prove in \cref{sec:RP  distribution of norm} is not surprising.
\end{remark*} 

\section{Random projection preserves distribution of norm}\label{sec:RP  distribution of norm}
Notice by central limit theorem, we have 
\[
\frac{\|X\|^2 - n}{\sqrt{n (\E x_1^4-1)}} \xrightarrow[]{n \to \infty} \cN (0,1)
\]

To understand if the random projected norm $\| SX\|^2$ has any type of convergence, it is necessary to find proper center and scale which are corresponding to first and second moments. Let us compute the mean first.
\begin{align*}
\E \| SX\|^2 -mn 
   = &  \E X^TS^T SX -mn \\
    = & \E \tr(XX^TS^TS) - mn  \\
    = & \tr(\E XX^T \E S^TS) - mn \\
    = & \tr(I_n\; m I_n) - mn \\
    =  & \;0
\end{align*}
For the variance, 
\begin{align*}
 & \E (\| SX\|^2  -mn)^2 \\
 = & \E (X^TS^T SX)^2 -m^2n^2 \\
 = & \E (X^TS^T SXX^TS^T SX) -m^2n^2 \\
    = & \E \left( \sum x_{i_1} S_{i_2,i_1} S_{i_2,i_3} x_{i_3} x_{j_1} S_{j_2,j_1} S_{j_2,j_3} x_{j_3} \right) - m^2n^2 
\end{align*}
The surviving terms must have even powers since first moments of the random variables are all 0. Therefore we only need to count four cases 
$\{i_1=i_3=j_1=j_3:=i\}$, $\{i_1=i_3:=i\neq j_1=j_3:=j\}$, $\{i_1=j_1:=i\neq i_3=j_3:=j\}$, $\{i_1=j_3:=i\neq i_3=j_1:=j\}$.\\
\begin{align*}
& \E \left( \sum x_{i_1} S_{i_2,i_1} S_{i_2,i_3} x_{i_3} x_{j_1} S_{j_2,j_1} S_{j_2,j_3} x_{j_3} \right)  \\
= & \E \sum_{i_2, j_2=1}^m \sum_{i=1}^n x_i S_{i_2,i} S_{i_2,i} x_{i} x_{i} S_{j_2,i} S_{j_2,i} x_{i}  
+  \E  \sum_{i_2, j_2=1}^m \sum_{\substack{i,j=1\\ i\ne j}}^n x_{i} S_{i_2,i} S_{i_2,i} x_{i} x_{j} S_{j_2,j} S_{j_2,j} x_{j}   \\
& + \E  \sum_{i_2, j_2=1}^m \sum_{\substack{i,j=1\\ i\ne j}}^n x_{i} S_{i_2,i} S_{i_2,j} x_{j} x_{i} S_{j_2,i} S_{j_2,j} x_{j}  
 + \E  \sum_{i_2, j_2=1}^m \sum_{\substack{i,j=1\\ i\ne j}}^n x_{i} S_{i_2,i} S_{i_2,j} x_{j} x_{j} S_{j_2,j} S_{j_2,i} x_{i} \\
= & \E \sum_{i_2, j_2=1}^m \sum_{i=1}^n x_i^4 S_{i_2,i}^2  S_{j_2,i}^2   
+ \E  \sum_{i_2, j_2=1}^m \sum_{\substack{i,j=1\\ i\ne j}}^n x_{i}^2 x_{j}^2 S_{i_2,i}^2 S_{j_2,j}^2 \\
& + \E  \sum_{k=1}^m \sum_{\substack{i,j=1\\ i\ne j}}^n x_{i}^2 x_{j}^2 S_{k,i}^2 S_{k,j}^2   
+ \E  \sum_{k=1}^m \sum_{\substack{i,j=1\\ i\ne j}}^n x_{i}^2 x_{j}^2  S_{k,i}^2 S_{k,j}^2  \\
& \text{where the last two terms dropped the zero terms $i_2\ne j_2$, } \\
= & (m^2n-mn+mn\E S_{11}^4)\E x_1^4 + m^2(n^2-n) + 2m(n^2-n) \\
= & m^2n^2 + (\E x_1^4-1) m^2n + 2mn^2 + mn[(\E S_{11}^4-1)\E x_1^4 -2] \\
= & m^2n^2 + \sigma^2 m^2n + 2mn^2 + \xi mn
\end{align*}
where $\sigma^2:= \E S_{11}^4-1, \; \xi:=[(\E S_{11}^4-1)\E x_1^4 -2]$. Therefore 
\begin{equation}\label{eqn: projected norm variance}
    \E \left[\left( \frac{X^TS^TSX - mn}{\sqrt{\sigma^2 m^2n+2mn^2+\xi mn}} \right)^2 \right] 
    =   1 
\end{equation}
Next we will show the centered and scaled projected norm actually also converges to a normal. That means distribution of the norm of a  vector (with independent entries) is also invariant under random projection.

\subsection{CLT for Random projection of norm}
\begin{theorem}
\label{thm:rand_proj norm CLT}
Given a random vector $X$ in $\R^n$ with i.i.d. entries
$$X=\mat[c]{x_1\\ \vdots\\ x_n}$$
Let $\E x_1=0, \E x_1^2=1, \E x_1^4 =1+\sigma^2 (0\le \sigma <\infty)$.
Consider a  random  matrix $S: \R^n \to \R^m$ with independent identically distributed entries $S_{i,j}$ with $ \E S_{i,j} = 0$,  $ \E S_{i,j}^2 = 1$ and $\E S_{1,1}^4<c<\infty$.
Further assume $S,X$ are all independent.
Define 
\[
    A(m,n): =\frac{\|SX\|^2 - mn}{\sqrt{\sigma^2 m^2n+2mn^2+ \xi mn}}
\]
where $\xi=[(\E S_{11}^4-1)\E x_1^4 -2]$. If $\frac{m}{n} \to 0$, then
\begin{align}
 A(m,n)\xrightarrow[]{m,n \to \infty} \cN(0,1) \label{eqn:rand_proj norm CLT}
\end{align}

\end{theorem}
\begin{remark*}
Before we proceed with the proof, it is worth mentioning that the random norm is a complicated sum of $mn^2$ correlated terms.
\[
X^TS^TSX = \sum_{k=1}^m\sum_{i,j=1}^n x_i S_{k,i} S_{k,j}x_j
\]
Therefore most analytical methods and tools fail to treat the quantity properly. For example,  characteristic function need independent property, Lindeberg swapping needs martingale property, and Stein's method needs exchangeable structure and precise control of first and second moments of conditional perturbed differences. So we are constrained to use a robust and universal approach, the moment method, which can be used to prove convergence to a limit law with known moments.

In the moment method we present below, we need to control the order of $m$ by $m \le o(n)$ because the counting procedure would be impossible to carry out if this is not the case. The scaling in this case is dominated by $\sqrt{2mn^2}$ which allows us to limit the significant terms in the moment calculation. However in simulations, we will see convergence to normal even when $m> n$. But we could not find a proof for the general case yet due to too many correlated terms. 
\end{remark*}

\begin{proof}
We will first note that a truncation argument will show it is sufficient to prove the same CLT result for bounded random variables. Details can be found in many standard moment method proof for CLT in many standard textbook (see for example \cite{tao2011topics} 2.2).

From now on, we assume all random variables are bounded, so that they have finite moments of all order which is very important in moment method. We will compute all moments of $A(m,n)$ in the limit and we expect all odd moments vanish and all even moments match with standard normal random variable.

The key idea is to separate the random norm into $m$ identically distributed but dependent random variables. Let $S_{k,.}$ be $k$-th row of $S$. Then $\|SX\|^2 = \sum_{k=1}^m (S_{k,.} X)^2$. Define 
\begin{align*}
 & L_k := \frac{(S_{k,.}X)^2 -n}{\sqrt{\sigma^2 mn+2n^2 + \xi n}} 
 =\frac{(\sum_{i,j=1}^n x_i S_{k,i} S_{k,j}x_j -n)}{\sqrt{\sigma^2 mn+2n^2 + \xi n}}  \\ 
\implies & A(m,n) = \frac{1}{\sqrt{m}} \sum_{k=1}^m L_k \\
\implies & \lim_{m,n\to \infty} \E A(m,n)^t
=\lim_{m,n\to \infty} \frac{1}{m^{t/2}} \E \left(\sum_{k=1}^m L_k \right)^t\\
\end{align*}

\vspace{2em}
Let us first record some moments properties of these identically distributed $L_k$.
\begin{enumerate}[label={(\arabic*)}]
    \item $\E L_k=0$.
    \item $\E L_k^2 = 1 -O(\frac{m}{n})$.
    \item For any fixed $ t\in \R $, there is a constant  $C_t$  independent of $m$ and $n$ so that
    \begin{equation}\label{eqn:L_k bounded moments}
        |\E L_k^t |\le C_t <\infty 
    \end{equation}
    
    \item The order of the expectation of a polynomial is determined by the number of singletons. 
    Given integer $q_1, \cdots q_r \ge 0$,
\begin{equation}
   \E[L_1^{q_1}\cdots L_r^{q_r} ] \le O\left(\frac{1}{\left({\sigma^2 m+2n}\right)^{d/2}}\right), 
   \qquad \text{where } d = \sum_{i=1}^r 1_{(q_i=1)}
   \label{eqn:Rand norm powers of 1 multiply poly}
\end{equation}
 Here $d$ is the total number of variables $L_i$ with multiplicity 1.
    \item For any fixed $r \in \R$, 
    \[ 
    \E (L_{1}^{2} \cdots L_{r}^{2})=\left(\frac{2n}{\sigma^2m + 2n} \right)^r + O(\frac{1}{n})
    \]
    which converges to 1 if $\frac{m}{n} \to 0$.
\end{enumerate}
We will prove one by one.
\begin{enumerate}[label={(\arabic*)}]
\item 
Obviously, $\E L_k = 0$. 

\item The variance $\E L_k^2= 1 -O(\frac{mn}{\sigma^2 mn+2n^2 + \xi n}) $  since 
\begin{align*}
    \E(\sum_{i,j=1}^n x_i S_{k,i} S_{k,j}x_j -n)^2 & = \E(\sum_{i,j=1}^n x_i S_{k,i} S_{k,j}x_j )^2 - n^2 \\
    & = \E(\sum_{i=1}^n x_i^2 S_{k,i}^2 )^2+ 2 \E[(\sum_{i=1}^n x_i^2 S_{k,i}^2 )(\sum_{\substack{i,j=1\\ i\ne j}}^n x_i S_{k,i} S_{k,j}x_j)]\\
    & \qquad +\E(\sum_{\substack{i,j=1\\ i\ne j}}^n x_i S_{k,i} S_{k,j}x_j)^2 - n^2 \\
    & = [n\E x_1^4 S_{1,1}^4 +n(n-1)] +0 + 2n(n-1) -n^2 \\
    & = 2 n^2 +O(n)
\end{align*}
Therefore  $\E L_k^2= 1-O(\frac{1}{n})$ when $m=o(n)$.
\item It is also true that $L_k$ has finite moments of all order which also hold in the limit.
We will use a careful counting procedure. First we notice
\begin{align}
 | \E L_k^t |& = \left| \E \left[ \frac{(\sum_{i,j=1}^n x_i S_{k,i} S_{k,j}x_j -n)}{\sqrt{\sigma^2 mn+2n^2 + \xi n}} \right]^t  \right| \nonumber \\
 & \le n^{-t} \left| \E \left[ \left(\sum_{i,j=1}^n x_i S_{k,i} S_{k,j}x_j -n \right)^t \right] \right|
 \label{eqn: Rand norm L_k^t finite}
\end{align}
Since $n = \sum_{i,j}\delta_{i,j}$ where $\delta_{i,j}=1$ when $i=j$ and $0$ otherwise. Then we expand the $t$-th moment on the right as polynomials.
\begin{align*}
    & \E \left[ \left(\sum_{i,j=1}^n x_i S_{k,i} S_{k,j}x_j -n \right)^t\right] \\
    = & \E \left[ \left(\sum_{i_1,j_1=1}^n (x_{i_1} S_{k,i_1} S_{k,j_1}x_{j_1} -\delta_{i_1,j_1}) \right)\cdots \left(\sum_{i_t,j_t=1}^n (x_{i_t} S_{k,i_t} S_{k,j_t}x_{j_t} -\delta_{i_t,j_t}) \right)\right] \\
    = & \sum_{i_1,j_1=1}^n \cdots \sum_{i_t,j_t=1}^n  \E [(x_{i_1} S_{k,i_1} S_{k,j_1}x_{j_1} -\delta_{i_1,j_1}) \cdots  (x_{i_t} S_{k,i_t} S_{k,j_t}x_{j_t} -\delta_{i_t,j_t}) ]
\end{align*}
Notice  $\E [(x_{i_1} S_{k,i_1} S_{k,j_1}x_{j_1} -\delta_{i_1,j_1}) \cdots  (x_{i_t} S_{k,i_t} S_{k,j_t}x_{j_t} -\delta_{i_t,j_t}) ] $ vanishes if cardinality of $\{i_1, j_1, \cdots i_t, j_t \}$ is greater than $t$ since that will have at least one singleton factor $\E S_{k,i}=0$. So the index set must collapse to a set of at most $t$ distinct indices. This means total number of nonzero terms is exactly $n^t$. As we are assume all variables are truncated to have bounded moments, we find
\[
\left| \E \left[ \left(\sum_{i,j=1}^n x_i S_{k,i} S_{k,j}x_j -n \right)^t\right] \right| \le 
O(n^t)
\]
Therefore plugging into \cref{eqn: Rand norm L_k^t finite} we find $|\E L_k^t |\le O(1)$.
This proves \cref{eqn:L_k bounded moments}.

\item \cref{eqn:Rand norm powers of 1 multiply poly} is an important property concerns the product of $L_k$'s. 
The order of the expectation of a polynomial $\E[L_1^{q_1}\cdots L_r^{q_r} ]$ is determined by the number of singletons $d = \sum_{i=1}^r 1_{(q_i=1)}$. 
More precisely, for any product of $L_k$ involving $d$ term of power 1, then the expected value is of order $(m+n)^{-d/2}$. In other words, each power 1 term contribute a factor of $(m+n)^{-1/2}$. we use an argument by conditioning  and careful counting. First, noticing $L_k$ conditioning on $X$ are independent,
\begin{align*}
    \E [L_1^{q_1} L_2^{q_2}\cdots L_r^{q_r} ]
    & = \E \left[(L_2^{q_2}\cdots L_r^{q_r})\E(L_1^{q_1}| X, S_{2,.},\cdots S_{r,.} ) \right] \\
    & = \E \left[(L_2^{q_2}\cdots L_r^{q_r})\E(L_1^{q_1}| X ) \right] \\
    & = \E \left[ \E \left((L_2^{q_2}\cdots L_r^{q_r})|X, S_{3,.},\cdots S_{r,.} \right) \E(L_1^{q_1}| X )  \right] \\
    & = \E \left[ (L_3^{q_3}\cdots L_r^{q_r})\E \left(L_2^{q_2}| X \right) \E(L_1^{q_1}| X )  \right] \\
    & \cdots \\
    & = \E \left[  \E(L_1^{q_1}| X ) \E \left(L_2^{q_2}| X \right) \cdots \E \left(L_r^{q_r}| X \right) \right] 
\end{align*}
Without loss of generality, assume the first $d$ variables are of multiplicity 1, namely 
\[
q_1 = q_2 \cdots = q_d =1
\]
To simplify notation, we denote  $\E (L_i | X) :=\mu_i, 1\le i\le d$. We would also only need the above conditioning argument up to $d$-th variable $L_d$. That is 
\begin{align*}
    \E [L_1\cdots L_d L_{d+1}^{q_{d+1}}\cdots L_r^{q_r} ]
    & = \E [\mu_1 \cdots \mu_d L_{d+1}^{q_{d+1}}\cdots L_r^{q_r} ]
\end{align*}
Then we apply Cauchy-Schwarz inequality, we find
\begin{align}
    \E [L_1\cdots L_d L_{d+1}^{q_{d+1}}\cdots L_r^{q_r} ]
    & \le \left( \E [\mu_1^2 \cdots \mu_d^2] \right)^{\frac{1}{2}}  \cdot \left( \E [L_{d+1}^{2q_{d+1}}\cdots L_r^{2q_r} ] \right)^{\frac{1}{2}} 
    \label{eqn:Rand norm powers of 1 Cauchy-Schwarz}
\end{align}

Then we notice the random variables $\mu_1, \cdots \mu_d$  (conditional expectation is also a random variable) are identical random variables (not just identically distributed).
\begin{align*}
  \mu_1 = \E \left( L_1| X\right)
  & = \E \left(  \frac{(\sum_{i,j=1}^n x_i S_{1,i} S_{1,j}x_j -n)}{\sqrt{\sigma^2 mn+2n^2 + \xi n}} \middle| X\right) \\
  & = \frac{1}{\sqrt{\sigma^2 mn+2n^2 + \xi n}}\left(  \sum_{i=1}^n x_i^2 \E S_{1,i}^2 +\sum_{\substack{i,j=1\\i\ne j}}^n x_ix_j \E(S_{1,i} S_{1,j}) -n \right) \\
  & = \frac{\sum_{i=1}^n x_i^2 -n}{\sqrt{\sigma^2 mn+2n^2 + \xi n}} \xrightarrow[]{\;a.s.\;} 0
\end{align*}
This shows $\mu_1$ does not depend on the index $1$ and indeed $\mu_1, \cdots, \mu_d$ are identical. As a side note, $\mu_1$ converges to 0 almost surely due to the strong law of large number. We would not need this fact though. Now we can analyze first term on the right hand side of \cref{eqn:Rand norm powers of 1 Cauchy-Schwarz},
\begin{align}
  \E \left[\mu_1^2 \cdots \mu_d^2 \right]
  & = \E \left[\mu_1^{2d}  \right] \nonumber \\
  & = \frac{\E \left(\sum_{i=1}^n x_i^2 -n\right)^{2d}}{ \left(\sigma^2 mn+2n^2 + \xi n \right)^{d} } \nonumber \\
  & = \frac{1}{ \left(\sigma^2 m+2n + \xi \right)^{d} }  \E \left(\frac{\sum_{i=1}^n x_i^2 -n}{\sqrt{n}}\right)^{2d} \nonumber \\
  &= O\left(\frac{1}{(\sigma^2 m+2n)^d}\right) \label{eqn: cauchy-schwarz L1 bound}
\end{align}
The last step, we used the the fact
\[
\E \left(\frac{\sum_{i=1}^n x_i^2 -n}{\sqrt{n}}\right)^{t} \le C_t <\infty
\]
that is due to independent sums $({\sum_{i=1}^n x_i^2 -n})$ in CLT has $t$-th moments of $O(n^{t/2})$ if $x_i$ has finite moments of all order (see \cite{brillinger1962note, von1965convergence}), which holds true due to its boundedness by truncation argument. 


On the other hand, applying Cauchy-Schwarz inequality repeatedly for the second term on the right hand side of \cref{eqn:Rand norm powers of 1 Cauchy-Schwarz}, we can bound it by a constant $c$ that does not depend on $m$ and $n$. Namely using the fact $L_k$ has bounded moments of all order (\cref{eqn:L_k bounded moments})
\begin{align}
    \E [L_{d+1}^{2q_{d+1}}\cdots L_r^{2q_r} ]
    & \le \left[ \E L_{d+1}^{4q_{d+1}} \E\left[ L_{d+2}^{4q_{d+2}}\cdots L_r^{4q_r}\right] \right]^{1/2} \nonumber \\
    & \le \left[\E L_{d+1}^{4q_{d+1}}\right]^{\frac{1}{2}} \left[\E L_{d+2}^{8q_{d+2}}\right]^{\frac{1}{4}} \cdots \left[\E L_r^{2^{r-d+1}q_r}\right]^{\frac{1}{2^{r-d}}} \nonumber\\
    & \le c <\infty \label{eqn:repeated cauchy schwartz show finite moments}
\end{align}
Therefore \cref{eqn:repeated cauchy schwartz show finite moments} combined with \cref{eqn:Rand norm powers of 1 Cauchy-Schwarz} and \cref{eqn: cauchy-schwarz L1 bound}, we obtain our desired result \cref{eqn:Rand norm powers of 1 multiply poly}.

\item 
To compute $\E (L_{1}^{2} \cdots L_{r}^{2}) $ we will use conditioning argument again.
\begin{align*}
    \E (L_{1}^{2} \cdots L_{r}^{2}) 
  & = \E \left[ \E (L_{1}^{2}|X) \cdots \E (L_{r}^{2}|X) \right] 
\end{align*}
Same as before, $\E (L_{1}^{2}|X),\cdots,\E (L_{r}^{2}|X)$ does not depend on the indices $1,\cdots, r$, and they are all identical random variables not just with same distribution. By definition,
\begin{align*}
  \E \left( L_1^2| X\right)
  =  \frac{\E  \left(  (\sum_{i,j=1}^n x_i S_{1,i} S_{1,j}x_j -n)^2 \middle| X\right) }{\sigma^2 mn+2n^2 + \xi n} 
\end{align*}
We shall simplify the numerator (denoted as $Q$),
\begin{align*}
  Q:= & \E   \left( (\sum_{i,j=1}^n x_i S_{1,i} S_{1,j}x_j -n)^2 \middle| X\right) \\
  = & \E \left( \left( \left(\sum_{i=1}^n (x_i^2 S_{1,i}^2 -1) \right) +\sum_{\substack{i,j=1\\i\ne j}}^n x_ix_j S_{1,i} S_{1,j} \right)^2 \middle| X \right)
\end{align*}
Expand the quadratic we find the cross terms vanish
\begin{align*}
    & 2  \E \left( \left(\sum_{k=1}^n (x_k^2 S_{1,k}^2 -1)\right) ( \sum_{\substack{i,j=1\\i\ne j}}^n x_ix_j S_{1,i} S_{1,j}) \middle| X \right)  \\
= & 2 \sum_{k=1}^n \sum_{\substack{i,j=1\\i\ne j}}^n \E \left[(x_k^2 S_{1,k}^2 -1)x_ix_j S_{1,i} S_{1,j} \middle| X \right] \\
= & 2 \sum_{k=1}^n \sum_{\substack{i,j=1\\i\ne j}}^n (x_k^2x_ix_j\E[S_{1,k}^2 S_{1,i} S_{1,j}]- x_ix_j\E[ S_{1,i} S_{1,j}] )
\end{align*}
 since each term $ \E [S_{1,k}^2 S_{1,i} S_{1,j}]=\E[ S_{1,i} S_{1,j}]=0$ as $i\ne j$. So we are left with
\begin{align*}
  Q = & \E \left( \left( \sum_{i=1}^n x_i^2 S_{1,i}^2 -n  \right)^2 \middle| X \right)
  + \E \left( \left( \sum_{\substack{i,j=1\\i\ne j}}^n x_ix_j S_{1,i} S_{1,j} \right)^2 \middle| X \right) \\
  = & \left( \sum_{i=1}^n\sum_{j=1}^n x_i^2x_j^2 \E S_{1,i}^2 S_{1,j}^2 - 2n \sum_{i=1}^n x_i^2\E S_{1,i}^2 + n^2 \right) \\
  & +  \left( \sum_{\substack{i,j=1\\i\ne j}}^n\sum_{\substack{i',j'=1\\i'\ne j'}}^n x_ix_jx_{i'}x_{j'} \E S_{1,i} S_{1,j} S_{1,i'} S_{1,j'}  \right)
\end{align*}
Notice the surviving terms in the second half are $\{i=i'\ne j=j' \}$ and $\{i=j'\ne j=i' \}$. Therefore, 
\begin{align*}
 Q = & \left(   \sum_{\substack{i,j=1\\i\ne j}}^n x_i^2x_j^2 +  \sum_{i=1}^n x_i^4\E S_{1,1}^4 - 2n \sum_{i=1}^n x_i^2 + n^2 \right) 
   +  \left( 2 \sum_{\substack{i,j=1\\i\ne j}}^n x_i^2 x_j^2  \right) \\
  = & 3\sum_{\substack{i,j=1\\i\ne j}}^n x_i^2x_j^2 - 2n \sum_{i=1}^n x_i^2 + n^2 + \sum_{i=1}^n x_i^4\E S_{1,1}^4 
\end{align*}
Since $\E(L_k^2|X)=Q/({\sigma^2 mn+2n^2 + \xi n} )$ for all $k$, we simplifies 
\begin{align*}
    \E (L_{1}^{2} \cdots L_{r}^{2}) 
  & = \E \left[ Q^r\right] /({\sigma^2 mn+2n^2 + \xi n} )^r
\end{align*}
To prove $\E (L_{1}^{2} \cdots L_{r}^{2}) \to 1$ it suffices to prove $\E \left[ Q / (2n^2)\right]^r \to 1 $ since $({\sigma^2 mn+2n^2 + \xi n} )^r$ is dominated by $(2n^2)^r$ as $m/n \to 0$.
\begin{align*}
    \frac{Q}{2n^2} & =  \frac{1}{{2n^2} }({3\sum_{\substack{i,j=1\\i\ne j}}^n x_i^2x_j^2 - 2n \sum_{i=1}^n x_i^2 + n^2 + \sum_{i=1}^n x_i^4\E S_{1,1}^4 } ) \\
& = \frac{1}{2}P +\frac{1}{2} + P_0
\end{align*}
where (to simplify notation) we denoted
\[
P := \left(\frac{1}{n^2}\sum_{i=1}^n x_i^2 (\sum_{j=1, j\ne i}^n  3x_j^2-2n)\right) ,\quad P_0 = \frac{1}{{2n^2}}\sum_{i=1}^n x_i^4\E S_{1,1}^4 
\]
First of all it is fairly easy to see $\E \left(\frac{Q}{2n^2} \right)^r<\infty$ for all $r\in \R$. This is because the total number of nonzero polynomial terms ($x_i^2, x_i^4, x_i^2x_j^2$) in $Q$ is $O(n^2)$ which will produce $O(n^{2r})$ polynomial terms for $Q^r$, and  $x_i$ has finite moment of all order. Now let us use induction to prove the moments are actually constant 1 in the limit. Suppose we have proved for all $k\le r$,
\[
\E \left(\frac{Q}{2n^2} \right)^k  \to 1, \quad \text{and } \E \left[ \left(\frac{Q}{2n^2} \right)^{k-1} P \right] \to 1
\]
It is easy to see this holds for the initial steps. Then we try to prove the next induction step namely 
\[
\E \left(\frac{Q}{2n^2} \right)^{r+1} \to 1, \quad \text{and } \E \left[ \left(\frac{Q}{2n^2} \right)^r P \right] \to 1
\]
Notice by linearity of expectation and Cauchy-Schwarz inequality,
\begin{align*}
    \E \left[ \left(\frac{Q}{2n^2} \right)^r P_0 \right] 
    &= \E \left[ \left(\frac{Q}{2n^2} \right)^r \left(  \frac{1}{{2n^2}}\sum_{i=1}^n x_i^4\E S_{1,1}^4  \right) \right] \\
    & = \frac{1}{{2n}}\E \left[ \left(\frac{Q}{2n^2} \right)^r  x_1^4\E S_{1,1}^4  \right] \\
    & \le \frac{1}{{2n}}\E \left[ \left(\frac{Q}{2n^2} \right)^{2r} \right]^{1/2}  \E \left[ x_1^8\E S_{1,1}^8  \right]^{1/2}  \to 0
\end{align*}
We find 
\begin{align*}
    \E \left(\frac{Q}{2n^2} \right)^{r+1} 
    & = \E \left[ \left(\frac{Q}{2n^2} \right)^r \left[ \frac{1}{2}P +\frac{1}{2} + P_0 \right] \right]  \\ & \to \frac{1}{2} \lim_{n\to \infty} \E \left[ \left(\frac{Q}{2n^2} \right)^r P \right] + \frac{1}{2} +0
\end{align*}
Now it suffices to show $\E \left[ \left(\frac{Q}{2n^2} \right)^r P \right] \to 1$. We expand again
\begin{align*}
    \E \left[ \left(\frac{Q}{2n^2} \right)^r P \right]
    & = \E \left[ \left(\frac{Q}{2n^2} \right)^{r-1} \left( \frac{1}{2}P^2 + \frac{1}{2}P +PP_0 \right) \right]    
\end{align*}
By induction hypothesis 
$\E \left[ \left(\frac{Q}{2n^2} \right)^{r-1}  \frac{1}{2}P  \right] = \frac{1}{2} $. And by Cauchy-Schwarz we know  $\E \left[ \left(\frac{Q}{2n^2} \right)^{r-1} PP_0  \right] \le O(\frac{1}{2n}) \to 0$. So it suffices to show $\E \left[ \left(\frac{Q}{2n^2} \right)^{r-1} P^2 \right] \to 1$. Repeat this argument $r$ times, we find it suffices to prove $\E P^{r+1} \to 1$.
And since this has to be true for every induction step, we indeed need to show 
\[
\E P^{k} \to 1, \quad \forall k\in \R
\]
Now we calculate
\begin{align*}
    \E P^k & = \E \left(\frac{1}{n^2}\sum_{i=1}^n  \sum_{j=1, j\ne i}^n  x_i^2(3x_j^2-\frac{2n}{n-1})\right)^k \\
    & = \frac{1}{n^{2k}} \sum_{i_1,\cdots, i_k =1}^n \sum_{j_1,\cdots, j_k =1, j_t\ne i_t}^n  \E \left( \prod_{t=1}^k x_{i_t}^2(3x_{j_t}^2-\frac{2n}{n-1})\right)
\end{align*}
The summations produce total $[n(n-1)]^k$ terms. Among them, there are total 
\[
n(n-1)\cdots (n-2k+1)= \prod_{s=0}^{2k-1}(n-s) = n^{2k} -O(n^{2k-1})
\]
terms that all indices are distinct, that is the cardinality $|\{i_1,\cdots, i_k, j_1,\cdots, j_k\}|=2k$. So that we can evaluate the expectation directly by independence. For $i_t \ne j_t, i_t\ne i_t', j_t\ne j_t'$
\begin{align*}
     \E \left( \prod_{t=1}^k x_{i_t}^2(3x_{j_t}^2-\frac{2n}{n-1})\right) 
     & = \prod_{t=1}^k \E x_{i_t}^2 \E (3x_{j_t}^2-\frac{2n}{n-1})  \\
     & =(1- \frac{2}{n-1})^k  = 1-O(\frac{1}{n}) \to 1
\end{align*}
For the remaining $n^k(n-1)^k - \prod_{s=0}^{2k}(n-s)=O(n^{2k-1})$ terms, the expectations of correlated variables are still bounded by some constant $c_k$. So we find 
\begin{align*}
    \E P^k 
    & = \frac{1}{n^{2k}} \prod_{s=0}^{2k}(n-s) + \frac{1}{n^{2k}} c_k\left( n^k(n-1)^k - \prod_{s=0}^{2k}(n-s)\right) \\
    &= 1-O(\frac{1}{n}) + \frac{1}{n^{2k}} c_k O(n^{2k-1}) \\
    & = 1+ O(\frac{1}{n})
    \to 1
\end{align*}
\end{enumerate}
\vspace{2em}
\bigskip
After establishing these properties, we can start the standard moment method for a CLT. By our definition of $L_k$, 
\begin{align*}
    \lim_{m,n\to \infty} \E A(m,n)^t
=\lim_{m,n\to \infty} \frac{1}{m^{t/2}} \E \left(\sum_{k=1}^m L_k \right)^t
\end{align*}
We first expand $\E \left(\sum_{k=1}^m L_k \right)^t$. Any term will have form $L_{p_1}^{q_1} \cdots L_{p_r}^{q_r}$ with the number of distinct indices $r: 1\le r \le t$ and positive integer powers satisfy $q_1+ \cdots +q_r = t$. If we group the terms by total number of distinct indices
\begin{align*}
    \E \left(\sum_{k=1}^m L_k \right)^t
    & = \sum_{r=1}^t \sum_{\substack{q_1+ \cdots +q_r = t\\ \{p_1,\cdots, p_r\}\subset \{1,\cdots,m\} }} c_{t,q_1,\cdots,q_r}\E (L_{p_1}^{q_1} \cdots L_{p_r}^{q_r})
\end{align*}
where $c_{t,q_1,\cdots,q_r}$ is the total number of orderings when we order $\{q_1$ number of index $p_1$, $\cdots$, $q_r$ number of index $p_r$, $q_1+ \cdots +q_r = t\}$ all together, which is 
\[
c_{t,q_1,\cdots,q_r} = \frac{t!}{q_1!\cdots q_r!}
\]
This constant only depend on $t$ and $q_1, \cdots, q_r$, and it may be upper bounded by $t^t$. 

Now we are going to analyze how much the terms contribute for each fixed $r$. In particular, we will show the only significant terms which will survive after scaling are when $r=\frac{t}{2}$ (if $t$ is odd, that means no surviving terms).

For any term $L_1^{q_1}\cdots L_r^{q_r} $ with $r> \frac{t}{2}$,  there are at least 
\begin{equation}
    2\left(r- \frac{t}{2} \right)= 2r-t \label{eqn:total singletons r > t/2}
\end{equation}
variables $L_i$ have multiplicity 1. It is true because increasing the length $r$ by 1 will create at least two singleton terms. For example if we want to increase length of $L_1^2L_2^2$ from 2 to 3, we would end up breaking a square term into two singletons so that we have $L_1^2L_2L_3$ or $L_1L_2^2L_3$. Formally, suppose that's not the case. Namely suppose there are only $s \le 2r-t -1$ variables $L_i$ of multiplicity 1. Then adding all the multiplicity we get 
\[
t=q_1+\cdots q_r \ge s+2(r-s) = 2r-s\ge t+1
\]
This is a contradiction.

Combining \cref{eqn:total singletons r > t/2} and \cref{eqn:Rand norm powers of 1 multiply poly}, each of the term when $r>\frac{t}{2}$ contribute at most 
\[
O\left(\frac{1}{\left({\sigma^2 m+2n}\right)^{d/2}}\right), 
   \qquad \text{where } d \ge 2r - t
\]
 For each fixed $r$, the total number of possible choices of $\{p_1,\cdots, p_r\}\subset \{1,\cdots,m\}$ is $\binom{m}{r} \le O(m^r)$. Then we also need to count total number of ways to generate $q_1+\cdots q_r=t$. That is we are looking at separating the integer $t$ into $r$ nonzero integers $q_1, \cdots, q_r$. Total number will be $\binom{t-1}{r-1} \le O(t^r)$ which only depend on $t$ (we can model it as separating $t$ stones into $r$ piles.  $\binom{t-1}{r-1}$ is due to the fact we can select $r-1$ separating positions out of $t-1$ spaces.).
 
 Therefore the total contribution for each fixed $r> \frac{t}{2}$ is 
 \begin{align*}
    &\quad  \frac{1}{m^{t/2}} \sum_{\substack{q_1+ \cdots +q_r = t\\ \{p_1,\cdots, p_r\}\subset \{1,\cdots,m\} }} c_{t,q_1,\cdots,q_r}\E (L_{p_1}^{q_1} \cdots L_{p_r}^{q_r}) \\
    & \le m^{-\frac{t}{2}}O(m^r) O(t^r) c_{t,q_1,\cdots,q_r} O\left(\frac{ 1 }{\left({\sigma^2 m+2n}\right)^{d/2}}\right) \\
    & \le O\left(m^{r-\frac{t}{2}} {\left({{\sigma^2m+ 2n}}\right)^{-d/2}}\right)  \\
    & \le O\left(m^{r-\frac{t}{2}}  \left({\sigma^2m+ 2n}\right)^{-r+\frac{t}{2}}\right)
\end{align*}
 Last step we used the fact $d\ge 2r-t$. Since we assumed $\frac{m}{n} \to 0$, we find the total contribution is bounded by
 \[
 O \left(\frac{m}{\sigma^2 m +2n}\right)^{r-\frac{t}{2}} \to 0, \quad \forall r > \frac{t}{2}
 \]
Therefore we can safely drop all cases of $r> \frac{t}{2}$. 
 
 For all cases $r<\frac{t}{2}$, each term $\E (L_{p_1}^{q_1} \cdots L_{p_r}^{q_r}) =O(1)$ depend only on $t$ by a repeated Cauchy-Schwarz argument similar to \cref{eqn:repeated cauchy schwartz show finite moments}. We find total contribution is
 \[
  m^{-\frac{t}{2}}O(m^r) O(t^r) c_{t,q_1,\cdots,q_r} O(1) \to 0, \quad \forall r < \frac{t}{2}
 \]
 This implies for $t$ odd
 \[
  \lim_{m,n\to \infty} \frac{1}{m^{t/2}} \E \left(\sum_{k=1}^m L_k \right)^t \to 0
 \]
 Now for even moments, the above analysis shows we only need to count $r=\frac{t}{2}$ (since contributions from $r<\frac{t}{2}$ and $r>\frac{t}{2}$ are both negligible). In this case, we are looking at separate the integer $t$ into $r=\frac{t}{2}$ positive integers $q_1, \cdots, q_r$. The total number will be $\binom{t-1}{r-1}$ which only depends on $t$ ($\binom{t-1}{r-1}$ is due to the fact we can select $r-1$ separating positions out of $t-1$ spaces). There are still many terms not significant. By \cref{eqn:Rand norm powers of 1 multiply poly}, any term has a $L_k$ of multiplicity 1 will contribute at most $O(n^{-\frac{1}{2}})$, thus we may drop these terms. Among these $\binom{t-1}{r-1}$  terms, there is only one term that every $L_k$ has multiplicity at least 2, which will survive the scaling namely
 \[
q_1= \cdots = q_r = 2 
\]
 In other words,
\begin{align*}
    \lim_{m,n\to \infty} \frac{1}{m^{t/2}} \E \left(\sum_{k=1}^m L_k \right)^t
    & = \lim_{m,n\to \infty} \frac{1}{m^{t/2}} \sum_{ \substack{r=t/2 \\ \{p_1,\cdots, p_r\}\subset \{1,\cdots,m\} } } c_{t,2,\cdots,2}\E (L_{p_1}^{2} \cdots L_{p_r}^{2}) \\
    & =\lim_{m,n\to \infty} \frac{1}{m^{t/2}}  \binom{m}{t/2}c_{t,2,\cdots,2} \E (L_{1}^{2} \cdots L_{t/2}^{2}) \\
    & = \lim_{m\to \infty} \frac{1}{m^{t/2}}  \binom{m}{t/2} \frac{t!}{2^{t/2}} \\
    & = \frac{t!}{\frac{t}{2}! 2^{t/2}}
\end{align*}
Which is exactly the even moments of standard normal random variable $\cN(0,1)$. One can see this by finding a recurrent relation between moments using moment generating function. $$\E[\cN(0,1)]^{2k+2} =(2k+1)\E[\cN(0,1)]^{2k}$$
\end{proof}

\subsection{Simulation}
We first give some simulations to show the random  projected norm converges to normal distribution.
\begin{figure}[!ht]
	\includegraphics[width=0.95\linewidth]{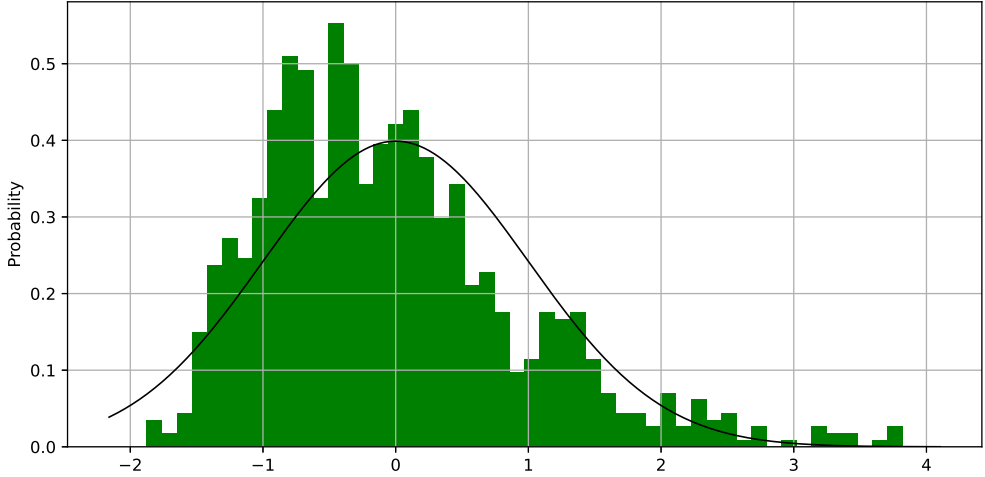}
	\caption{Random projected norm (m=10, n=100)}
	\label{fig: norm_1 Rand Proj}
\end{figure}
\begin{figure}[!ht]
    \includegraphics[width=0.95\linewidth]{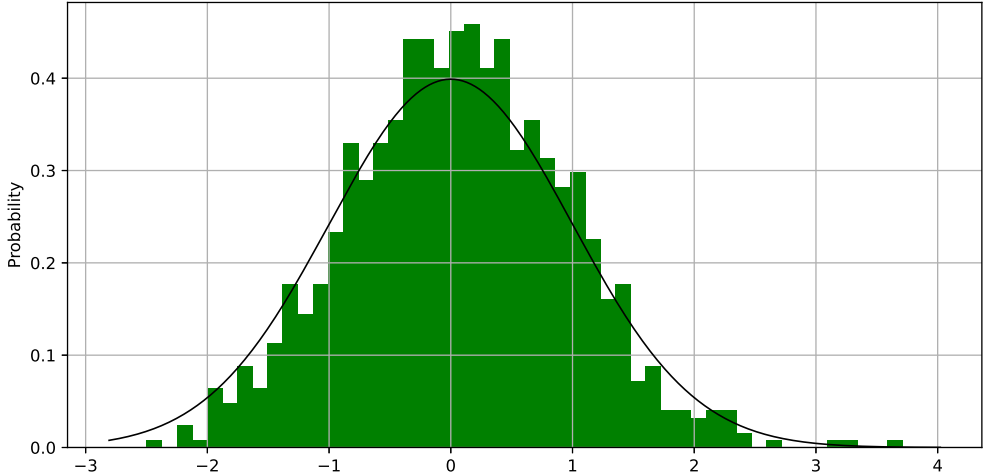}
    \caption{Random projected norm (m=500, n=5000)}
    \label{fig: norm_2 Rand Proj}
\end{figure}

\cref{fig: norm_1 Rand Proj} and \cref{fig: norm_2 Rand Proj} plotted histograms of 1000 samples of the projected norm
\[
\frac{1}{\sqrt{\sigma^2 m^2n+2mn^2+ \xi mn}} X^TS^TSX
\]
with different dimension settings. The random variables we used for entries of $X, S$ are standard normal random variables. As dimension $m, n$ increases, the convergence improves.

Next we give simulations for random embedded norms where $m> n$. Even though we do not have a CLT in this setting but the Bernstein type (mixed sub-Gaussian and sub-exponential) concentration behavior we proved in \cref{sec:concen norm by JL} is still relevant.
\begin{figure}[!ht]
	\includegraphics[width=0.95\linewidth]{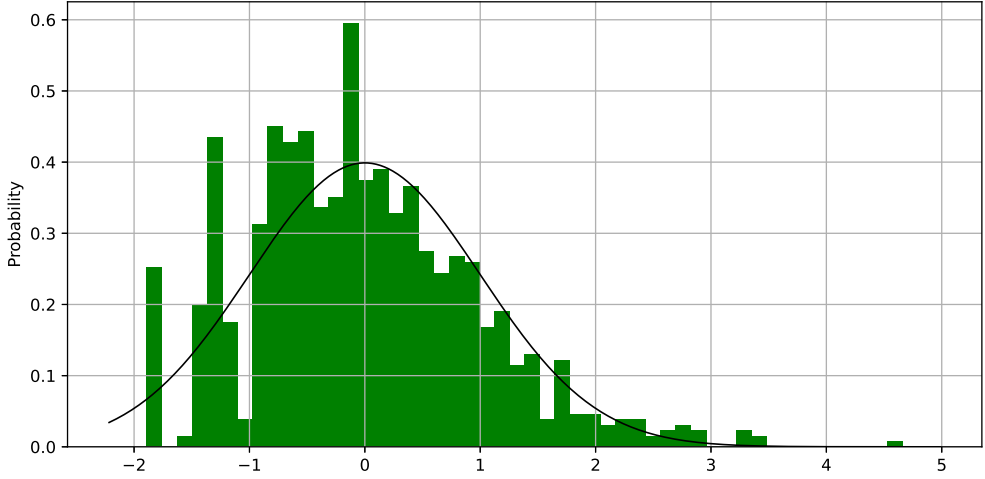}
	\caption{Random embedded norm (m=200, n=20)}
	\label{fig: norm_1E Rand Proj}
\end{figure}
\begin{figure}[!ht]
    \includegraphics[width=0.95\linewidth]{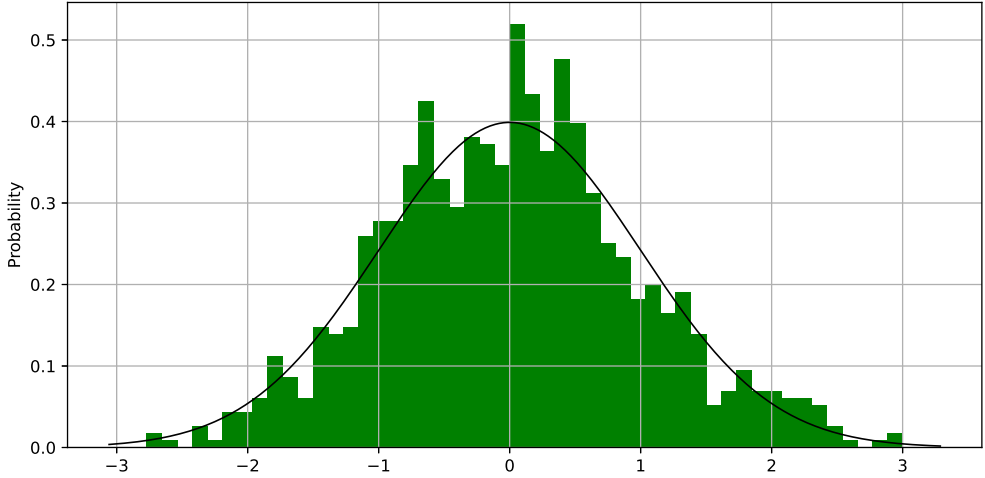}
    \caption{Random embedded norm (m=2000, n=200)}
    \label{fig: norm_2E Rand Proj}
\end{figure}

\cref{fig: norm_1E Rand Proj} and \cref{fig: norm_2E Rand Proj} 
plotted histograms of 1000 samples of the random embedded norm. The random variables we used for $X, S, Z$  take discrete values $\{-2.5, 0, 2.5\}$ with probability $\{0.08, 0.84, 0.08\}$. The kurtosis of such random variable is $6.25$ which is larger than standard normal random variable. Again as dimension $m, n$ increases, the histogram converges to a standard normal shape.


\subsection{Possibility of extending CLT for random projection to random embedding}
\cref{fig: norm_1E Rand Proj} and \cref{fig: norm_2E Rand Proj} shows it is very likely there is a CLT for random embedded norms where $m\ge O(n)$. Of course the moment computation would fail because of the complicated dependence structure.

If we view the random transformed norm $\|SX\|^2$ as a trace function on the spectral of product of random matrices, there are potential ways from random matrix theory to overcome the difficult of too many correlated random variables when $m\ge O(n)$.
\[
X^TS^TSX = \tr(X^TS^TSX) =\tr(S^TSXX^T) 
\]
It is clear the random matrix $\frac{1}{\sqrt{n}} XX^T$ has one nonzero eigenvalue which is $\lambda_1 = \frac{1}{\sqrt{n}}\sum_{i=1}^n x_i^2 = O(\sqrt{n})$. 

From random matrix analysis \cite{marvcenko1967distribution}, it is known that the empirical spectral distribution of $A_n:=\frac{1}{m}S^TS$ converges to the celebrated Mar{\v c}enko-Pastur law depends on the parameter $\frac{m}{n} \to c$, $c \in (0,\infty)$.   
\[
 \frac{1}{n} \sum_i \delta_{\lambda_i(A_m)}  \xrightarrow{m\to \infty}  \mu_{M,c}
\]
where $\lambda_i(A_n)$ are eigenvalues of $A_n$ and $\delta_{\lambda_i(A_n)}$ is the Dirac delta function, $\mu_{M,c}$ is the Mar{\v c}enko-Pastur probability measure. Moreover, \cite{yin1986limiting, silverstein1995strong, bai2010spectral} showed if a non-negative definite $n\times n$ random matrix $B_n$ has a deterministic limiting distribution $F^B$, then one can characterize the limiting spectral distribution of the product, $A_n B_n$, converges in distribution to probability distribution $F$.
\[
 \frac{1}{m} \sum_i \delta_{\lambda_i(A_nB_n)}  \xrightarrow{n\to \infty}  \mu_{F}
\]
One may start thinking if it is possible to apply the result of product of random matrices to our problem. Obviously, one would replace $B_n$ with  $\frac{1}{\sqrt{n}}XX^T$. However the spectral distribution of $B_n$ does not converge properly since $\lambda_1 = O(\sqrt{n})$. And our CLT result is actually on another level of details. One has to first center and standardize spectral of $B_n$, then see how the fluctuation is interacting with the spectral of $A_n$. In fact  $A_nB_n$ has only one nonzero eigenvalue, we are actually looking at distribution of this single eigenvalue, which usually requires very different techniques to compute. The extreme eigenvalues of full rank random matrices usually converges to Tracy-Widom distribution \cite{tracy1994level, johnstone2001distribution}.  In our case, we are looking at a version of this type but the random matrix has certain structure of rank one.


\section{Discussion on rate of convergence }\label{sec:Open questions}
In this section we discuss the rate of convergence for the projected or embedded norm. Based on some detailed calculation, we believe the following conjecture should be true.

\begin{conjecture}[Random projection of norm rate of invariance] \label{conj:rand_proj norm}
Given a random vector $X$ in $\R^n$ with i.i.d. entries
$$X=\mat[c]{x_1\\ \vdots\\ x_n}$$
Let $\E x_1=0, \E x_1^2=1, \E x_1^4 =1+\sigma^2 (0< \sigma <\infty)$.
Consider a  random  matrix $S: \R^n \to \R^m$ with independent entries  and $ \E S_{i,j} = 0$ and  $ \E S_{i,j}^2 = 1$.
Further assume $S,X$ are all independent and $\E S_{1,1}^8\vee \E |x_1|^6<c<\infty$. Also let $G$ be a standard normal random variable. Then we have
\begin{align}
    \sup_t \left| \p\left( \frac{\left(X^TS^TSX - mn\right)}{\sqrt{\sigma^2 m^2n+2mn^2}}  < t \right) - \p\left( G < t \right) \right| 
    \le & O\left(\frac{1}{\sqrt{n}}+ \frac{1}{\sqrt{m}} \right)  \label{eqn:rand_proj norm to normal} \\
\sup_t  \left| \p\left( \frac{\left(X^TS^TSX - mn\right)}{\sqrt{\sigma^2 m^2n+2mn^2}}  < t \right) - \p\left( \frac{( X^T X  - n)}{\sigma\sqrt{n}}   < t \right) \right| 
    \le & O\left(\frac{1}{\sqrt{n}}+ \frac{1}{\sqrt{m}} \right)  \label{eqn:rand_proj norm invariance}
\end{align}
\end{conjecture}
\begin{remark*}
If we use Berry-Essen theorem for random sequence $X^2-1 = [x_1^2-1,\cdots, x_n^2-1]$ (the assumptions in Berry-Essen are satisfied since $x_i^2 -1$ are i.i.d., $\E x_1^2-1 = 0$, $\E(x_1^2 -1)^2=\sigma^2$ and  $\E|x_1^2 -1|^3<\infty$.) we find 
\[
\sup_t  \left| \p\left( G < t \right) - \p\left( \frac{( X^T X  - n)}{\sigma\sqrt{n}}   < t \right) \right| 
    \le  O\left(\frac{1}{\sqrt{n}} \right) 
\]
Then it is tempting to use techniques similar with   \cite{duan2021invariance} to prove \cref{eqn:rand_proj norm to normal}, then by triangle inequality one concludes \cref{eqn:rand_proj norm invariance}. However, the techniques in \cite{duan2021invariance} heavily relies on the fact that we can separate the quantity of interests into two independent parts $X$  and $S^TSZ$. For this conjecture on the rate of norm invariance, there is no such luxury property that we can exploit.  
\end{remark*}


There are some examples suggesting this rate is the correct order. We will try to analyze in detail to see how much distortion is introduced in the projected norm with two examples. From the variance calculation \cref{eqn: projected norm variance}, we know there is an error term at least the order $O(\frac{1}{\sqrt{m+n}})$. Then let us analyze two special cases $m=1, n\to \infty$ and $m\to \infty, n=1$.

For $m=1, n\to \infty$, we find 
\begin{align*}
    \frac{\left(X^TS^TSX - mn\right)}{ \sqrt{ \sigma^2 m^2n+2mn^2}}   
    & = \frac{\left(\sum_{i=1}^n S_{1,i} x_i\right)^2 - n}{ \sqrt{m}\sqrt{ \sigma^2n +2n^2}}  \\
    & = \frac{1}{\sqrt{m}}\sqrt{\frac{n^2}{\sigma^2n +2n^2}}\left[ \left(\frac{\sum_{i=1}^n S_{1,i} x_i}{\sqrt{n}}\right)^2 -1 \right] \\
    & \approx \frac{1}{\sqrt{m}} \sqrt{\frac{n^2}{\sigma^2n +2n^2}}\left[ \cN(0,1)^2 -1 \right] \\
    & = \frac{1}{\sqrt{2m}}[\cN(0,1)^2 -1]
\end{align*}
Since $m=1$, $\cN(0,1)^2 -1$ differ from $\cN(0,1)$ by $O(1)$, we see the error term is $O(\frac{1}{\sqrt{m}})$. 

For $m\to \infty, n=1$, similarly we compute 
\begin{align*}
    \frac{\left(X^TS^TSX - mn\right)}{ \sqrt{ \sigma^2 m^2n+2mn^2}}   
    & = \frac{x_1^2 \left(\sum_{i=1}^m S_{i,1}^2 \right) - m}{ \sqrt{n}\sqrt{ \sigma^2 m^2 +2m}}  \\
    & = \frac{x_1^2 \left(\sum_{i=1}^m S_{i,1}^2-m \right) +m(x_1^2-1)}{ \sqrt{n}\sqrt{ \sigma^2 m^2 +2m}}  \\
    & \approx \frac{x_1^2}{\sqrt{mn}} \cN(0,1) + \frac{x_1^2-1}{\sqrt{\sigma^2 n}} \\
    & \approx 0 + \frac{x_1^2-1}{\sqrt{\sigma^2 n}}
\end{align*}
In this case the error term is on the scale of $O(\frac{1}{\sqrt{n}})$ since $x_1^2-1$ differs from $\cN(0,1)$ by O(1). For large $m$ and $n$, we believe both $O(\frac{1}{\sqrt{m}})$ and $O(\frac{1}{\sqrt{n}})$ are necessary. 

\bigskip

\vskip 0.2in
\bibliography{Invariancenorm}
\bibliographystyle{plain}
\end{document}